\documentclass[12pt]{amsart}
\usepackage {amsfonts,amscd,amssymb,amsmath,tabularx}
\usepackage{enumerate}

\newtheorem {prop*}{Proposition}
\newtheorem {theoreme*}[prop*]{Th\'eor\`eme}

\newtheorem {prop**}{Proposition}
\newtheorem {thm**}[prop**]{Theorem}

\newcommand\zset{\mathbb Z}

\newcommand\ilie {\mathfrak {i}}

\newcommand\glie {\mathfrak {g}}
\newcommand\hlie {\mathfrak {h}}
\newcommand\plie {\mathfrak {p}}

\newcommand \calf {{\mathcal F}}

\newcommand \cala {{\mathcal A}}

\newcommand \caln {{\mathcal N}}

\def\vtrait{\vrule height 10pt depth 6pt} 
\def\sboite#1{\vtrait\hbox to 16pt{\hfil #1\hfil}} 
\def\hboite#1{\hbox to 16pt{\hfil #1\hfil}} 

\def\vt{\vrule height 10pt depth 6pt} 
\def\sb#1{\vt\hbox to 36pt{\hfil #1\hfil}} 
\def\hb#1{\hbox to 36pt{\hfil #1\hfil}}

\newenvironment{Dynkin}{\setlength{\unitlength}{1.5pt}\begin{array}{l}}{\end{array}}

\newcommand{\Dbloc}[1]{\begin{picture}(20,20)#1\end{picture}}

\newcommand{\Dcirc}{\put(10,10){\circle{4}}}

\newcommand{\Deast}{\put(12,10){\line(1,0){8}}}
\newcommand{\Dwest}{\put(8,10){\line(-1,0){8}}}

\newcommand{\Ddoubleeast}{\put(10,12){\line(1,0){10}}\put(10,8){\line(1,0){10}}}
\newcommand{\Ddoublewest}{\put(10,12){\line(-1,0){10}}\put(10,8){\line(-1,0){10}}}

\newcommand{\Dtext}[2]{\makebox(20,20)[#1]{\scriptsize $#2$}}

\newcommand{\Dleftarrow}{\hskip-5pt{\makebox(20,20)[l]{\Large$<$}}\hskip-25pt}
\newcommand{\Drightarrow}{\hskip-5pt{\makebox(20,20)[l]{\Large$>$}}\hskip-25pt}
\begin{document}

\title[Enumeration for exceptional types]{Enumeration of ad-nilpotent ideals of parabolic subalgebras for exceptional types}
\author{C\'eline RIGHI}
\address{UMR 6086 CNRS, D\'epartement de Math\'ematiques, Universit\'e de Poitiers, T\'el\'eport 2 - BP
  30179, Boulevard Marie et Pierre Curie, 86962 Futuroscope
  Chasseneuil Cedex, France}
\email{celine.righi@math.univ-poitiers.fr}

\begin{abstract}
Using GAP4, we determine the number of ad-nilpotent and abelian ideals of a parabolic subalgebra of a simple Lie algebra of exceptional types $E$, $F$ or $G$. 
\end{abstract}

\maketitle
\thispagestyle{empty}

\section{Introduction}
Let $\glie$ be a complex simple Lie algebra of rank $l$. Let $\hlie$ be a  Cartan 
subalgebra of $\glie$ and $\Delta$ the associated root system. We fix a system of
positive roots $\Delta^+$. Denote by $\Pi=\{\alpha_1,
\dots, \alpha_l\}$ the corresponding set of simple roots. For each 
$\alpha \in \Delta$, let $ \mathfrak{g}_{ \alpha}$ 
be the root space of $\glie$ relative to $\alpha$.

For  $I \subset \Pi$, set $\Delta_I=\zset I \cap\Delta$.
We fix the corresponding standard parabolic subalgebra :   
$$
\plie_I=\hlie\oplus\left(\bigoplus\limits_{\alpha \in \Delta_I \cup \Delta^+}\glie_{\alpha}\right).
$$
An ideal $\ilie $ of $\plie_I$ is ad-nilpotent if and only if for all $x\in \ilie$, $ad_{\plie_I} x$ is nilpotent. Since any ideal of $\plie_I$ is $\hlie$-stable, we can deduce easily that an ideal is ad-nilpotent if and only if it is nilpotent. Moreover, we have $\ilie =\bigoplus\limits_{\alpha \in \Phi} \glie_{\alpha}$, for some subset $\Phi\subset \Delta^+ \setminus \Delta_I$. 

In this paper, we determine the number of ad-nilpotent and abelian ideals of a parabolic subalgebra when $\glie$ is of exceptional types $E$, $F$ and $G$. This is done by using GAP4. First, we recall some results of \cite{these} and \cite{R} on ad-nilpotent ideals. Then, we explain how we used these results for the programming and we give the tables of the enumeration in each exceptional type.

\section{Enumeration of ad-nilpotent ideals}\label{section_antichaine}

Let $I \subset \Pi$ and $\ilie$ be an ad-nilpotent ideal of $\plie_I$. We set 
$$
\Phi_{\ilie} =  \{ \alpha \in \Delta^+ \setminus \Delta_I;\ 
\mathfrak{g}_{\alpha } \subseteq \ilie \}.
$$ 
Then $\ilie = 
\bigoplus_{\alpha \in \Phi_{\ilie} } \mathfrak{g}_{\alpha}$ and if 
$\alpha \in \Phi_{\ilie}$, $\beta \in \Delta^+\cup \Delta_I$ are such that $\alpha +\beta
\in \Delta^+$, then $\alpha +\beta \in \Phi_{\ilie}$. 

Conversely, set 
$$ 
\calf_I = \{ \Phi \subset \Delta^+\setminus \Delta_I;\mbox{if } \alpha \in \Phi,
\beta \in \Delta^+\cup \Delta_I, \alpha +\beta \in \Delta^+, 
\mbox{then} \ \alpha +\beta\in \Phi \}.
$$
Then for $\Phi \in \calf_I$, $\ilie_{\Phi} = \bigoplus_{\alpha \in
\Phi} \mathfrak{g}_{\alpha}$ is an ad-nilpotent ideal of $\plie_I$. 

We obtain therefore a bijection
$$
\begin{array}{rcl}
\{\mbox{ad-nilpotent ideals of } \plie_I \}&  \rightarrow&  \calf_I ,\\ 
\ilie  &\mapsto&  \Phi_{\ilie}.
\end{array}
$$

Recall the following partial order on $\Delta^+$ : 
$\alpha \leqslant \beta$ if $\beta -\alpha$ is a 
sum of positive roots. Then it is easy to see that $\Phi \in
\calf_{\emptyset}$ 
if and only if for all $\alpha \in \Phi, \beta \in \Delta^+$, such
that $\alpha \leqslant \beta$, we have $\beta \in \Phi$. We can now define, for $\Phi\in \calf_{\emptyset}$ : 
$$
\Phi_{min}=\{\beta\in\Phi;
\beta-\alpha\not\in \Phi, \mbox{ for all }\alpha\in \Delta^+\}
$$
the set of minimal roots of $\Phi$ for $\leqslant$.

A set of positive roots $A$ is called an antichain of 
$(\Delta^+,\leqslant)$ if all the roots in $A$ are pairwise non comparable with respect to $\leqslant$. It is clear that $\Phi_{min}$ is an antichain. Conversely, let $A$ be an antichain. Set $\Phi=\{\beta\in\Delta^+; \mbox{ there exists } \alpha\in A \mbox{ such that } \alpha\leqslant \beta \}$, then $\Phi\in\calf_{\emptyset}.$ So each antichain corresponds to an element of $\calf_{\emptyset}$.

\bigskip
To obtain the enumeration of ad-nilpotent ideals for exceptional types, we first determine the set of antichains in $\Delta^+$. Then, by the considerations above, we obtain from the antichains all the elements of $\calf_{\emptyset}$. 

Next, to check if an element $\Phi\in\calf_{\emptyset}$ is an element of $\Phi\in\calf_{\{\alpha\}}$, for $\alpha\in\Pi$, it is enough to check that $\{\beta-\alpha; \alpha\in\Phi\}\cap (\Delta^+\cup\{0\}) \subset\Phi$. Thus, for $I\subset\Pi$, we obtain that the elements of $\calf_I$ are those of $\calf_{\emptyset}$ satisfying the condition above for each $\alpha\in I$. 

Let $\Phi\in\calf_I$, then $\Phi$ corresponds to an abelian ideal if and only if $\Phi^2=\emptyset$. Since the roots in $\Phi^2$ corresponds to those which are in the derived subalgebra of the ideal which corresponds to $\Phi$, it is also an ideal. So $\Phi^2\not=\emptyset$ if and only if the highest root $\theta$ is an element of $\Phi^2$. Then, to check if $\Phi$ is abelian, we need only to check that $\theta\not\in\Phi^2$.

\bigskip
Let $\caln_I$ be the set of ad-nilpotent ideals of $\plie_I$ and let $\cala b_I$ be the set of abelian ideals of $\plie_I$. The following tables give the cardinality of $\caln_I$ and $\cala b_I$ for each type $E$, $F$ and $G$. The subset $I$ of $\Pi$ is described by the symbol $\bullet$ in the Dynkin diagram without arrow which corresponds to the type we consider. 

For example, if we consider $E_6$, the diagram which corresponds to $I=\{\alpha_2,\alpha_5\}$ is 
$$

$$
where we use the numberings of \cite{TY}.

\medskip
The orientations for the diagram in type $F$ and $G$ are those in \cite{TY}.
\newpage
%% Table pour E6

\begin{center}
Table $E_{6}$
\end{center}

$$
\begin{array}{||c|c|c||c|c|c||}
\hline
I & \sharp\mathcal{N}_{I} & \sharp \mathcal{A}b_{I}  & 
I &  \sharp\mathcal{N}_{I} &  \sharp\mathcal{A}b_{I} \\
\hline
\begin{smallmatrix}& & \circ\\
\circ &\circ &\circ &\circ &\circ &\\
\end{smallmatrix}
&833&64& 
\begin{smallmatrix}& & \circ\\
\bullet &\circ &\circ &\circ &\circ &\\
\end{smallmatrix}
&197&21\\
\hline
\begin{smallmatrix}& & \bullet\\
\circ &\circ &\circ &\circ &\circ &\\
\end{smallmatrix}
&201&40& 
\begin{smallmatrix}& & \circ\\
\circ &\bullet &\circ &\circ &\circ &\\
\end{smallmatrix}
&255&40\\
\hline
\begin{smallmatrix}& & \circ\\
\circ &\circ &\bullet &\circ &\circ &\\
\end{smallmatrix}
&323&41& 
\begin{smallmatrix}& & \circ\\
\circ &\circ &\circ &\bullet &\circ &\\
\end{smallmatrix}
&255&40\\
\hline
\begin{smallmatrix}& & \circ\\
\circ &\circ &\circ &\circ &\bullet &\\
\end{smallmatrix}
&197&21& 
\begin{smallmatrix}& & \bullet\\
\bullet &\circ &\circ &\circ &\circ &\\
\end{smallmatrix}
&60&16\\
\hline
\begin{smallmatrix}& & \circ\\
\bullet &\bullet &\circ &\circ &\circ &\\
\end{smallmatrix}
&51&10& 
\begin{smallmatrix}& & \circ\\
\bullet &\circ &\bullet &\circ &\circ &\\
\end{smallmatrix}
&81&16\\
\hline
\begin{smallmatrix}& & \circ\\
\bullet &\circ &\circ &\bullet &\circ &\\
\end{smallmatrix}
&68&16& 
\begin{smallmatrix}& & \circ\\
\bullet &\circ &\circ &\circ &\bullet &\\
\end{smallmatrix}
&56&8\\
\hline
\begin{smallmatrix}& & \bullet\\
\circ &\bullet &\circ &\circ &\circ &\\
\end{smallmatrix}
&82&24& 
\begin{smallmatrix}& & \bullet\\
\circ &\circ &\bullet &\circ &\circ &\\
\end{smallmatrix}
&60&22\\
\hline
\end{array}
$$

\newpage
%% Table E7

\begin{center}
Table $E_{7}$
\end{center}

$$
\begin{array}{||c|c|c||c|c|c||}
\hline
I & \sharp\mathcal{N}_{I} & \sharp \mathcal{A}b_{I}  & 
I &  \sharp\mathcal{N}_{I} &  \sharp\mathcal{A}b_{I} \\
\hline
\begin{smallmatrix}& & \circ\\
\circ &\circ &\circ &\circ &\circ &\circ &\\
\end{smallmatrix}
&4160&128& 
\begin{smallmatrix}& & \circ\\
\bullet &\circ &\circ &\circ &\circ &\circ &\\
\end{smallmatrix}
&837&70\\
\hline
\begin{smallmatrix}& & \bullet\\
\circ &\circ &\circ &\circ &\circ &\circ &\\
\end{smallmatrix}
&980&78& 
\begin{smallmatrix}& & \circ\\
\circ &\bullet &\circ &\circ &\circ &\circ &\\
\end{smallmatrix}
&1251&73\\
\hline
\begin{smallmatrix}& & \circ\\
\circ &\circ &\bullet &\circ &\circ &\circ &\\
\end{smallmatrix}
&1600&84& 
\begin{smallmatrix}& & \circ\\
\circ &\circ &\circ &\bullet &\circ &\circ &\\
\end{smallmatrix}
&1385&73\\
\hline
\begin{smallmatrix}& & \circ\\
\circ &\circ &\circ &\circ &\bullet &\circ &\\
\end{smallmatrix}
&1076&70& 
\begin{smallmatrix}& & \circ\\
\circ &\circ &\circ &\circ &\circ &\bullet &\\
\end{smallmatrix}
&879&32\\
\hline
\begin{smallmatrix}& & \bullet\\
\bullet &\circ &\circ &\circ &\circ &\circ &\\
\end{smallmatrix}
&261&42& 
\begin{smallmatrix}& & \circ\\
\bullet &\bullet &\circ &\circ &\circ &\circ &\\
\end{smallmatrix}
&202&35\\
\hline
\begin{smallmatrix}& & \circ\\
\bullet &\circ &\bullet &\circ &\circ &\circ &\\
\end{smallmatrix}
&358&47& 
\begin{smallmatrix}& & \circ\\
\bullet &\circ &\circ &\bullet &\circ &\circ &\\
\end{smallmatrix}
&314&45\\
\hline
\begin{smallmatrix}& & \circ\\
\bullet &\circ &\circ &\circ &\bullet &\circ &\\
\end{smallmatrix}
&261&42& 
\begin{smallmatrix}& & \circ\\
\bullet &\circ &\circ &\circ &\circ &\bullet &\\
\end{smallmatrix}
&215&24\\
\hline
\end{array}
$$

\newpage

\begin{center}
Table $E_{7}$ (cont'd)
\end{center}

$$
\begin{array}{||c|c|c||c|c|c||}
\hline
I & \sharp\mathcal{N}_{I} & \sharp \mathcal{A}b_{I}  & 
I &  \sharp\mathcal{N}_{I} &  \sharp\mathcal{A}b_{I} \\
\hline
\begin{smallmatrix}& & \bullet\\
\bullet &\bullet &\bullet &\circ &\circ &\circ &\\
\end{smallmatrix}
&21&8& 
\begin{smallmatrix}& & \bullet\\
\bullet &\bullet &\circ &\bullet &\circ &\circ &\\
\end{smallmatrix}
&39&13\\
\hline
\begin{smallmatrix}& & \bullet\\
\bullet &\bullet &\circ &\circ &\bullet &\circ &\\
\end{smallmatrix}
&34&13& 
\begin{smallmatrix}& & \bullet\\
\bullet &\bullet &\circ &\circ &\circ &\bullet &\\
\end{smallmatrix}
&26&8\\
\hline
\begin{smallmatrix}& & \bullet\\
\bullet &\circ &\bullet &\bullet &\circ &\circ &\\
\end{smallmatrix}
&34&12& 
\begin{smallmatrix}& & \bullet\\
\bullet &\circ &\bullet &\circ &\bullet &\circ &\\
\end{smallmatrix}
&41&16\\
\hline
\begin{smallmatrix}& & \bullet\\
\bullet &\circ &\bullet &\circ &\circ &\bullet &\\
\end{smallmatrix}
&30&11& 
\begin{smallmatrix}& & \bullet\\
\bullet &\circ &\circ &\bullet &\bullet &\circ &\\
\end{smallmatrix}
&34&13\\
\hline
\begin{smallmatrix}& & \bullet\\
\bullet &\circ &\circ &\bullet &\circ &\bullet &\\
\end{smallmatrix}
&35&12& 
\begin{smallmatrix}& & \bullet\\
\bullet &\circ &\circ &\circ &\bullet &\bullet &\\
\end{smallmatrix}
&23&8\\
\hline
\begin{smallmatrix}& & \circ\\
\bullet &\bullet &\bullet &\bullet &\circ &\circ &\\
\end{smallmatrix}
&23&11& 
\begin{smallmatrix}& & \circ\\
\bullet &\bullet &\bullet &\circ &\bullet &\circ &\\
\end{smallmatrix}
&29&14\\
\hline
\begin{smallmatrix}& & \circ\\
\bullet &\bullet &\bullet &\circ &\circ &\bullet &\\
\end{smallmatrix}
&23&10& 
\begin{smallmatrix}& & \circ\\
\bullet &\bullet &\circ &\bullet &\bullet &\circ &\\
\end{smallmatrix}
&33&14\\
\hline
\end{array}
$$

\newpage
%% Tables pour E8

\begin{center}
Table $E_{8}$
\end{center}

$$
\begin{array}{||c|c|c||c|c|c||}
\hline
I & \sharp\mathcal{N}_{I} & \sharp \mathcal{A}b_{I} &
I & \sharp\mathcal{N}_{I} & \sharp \mathcal{A}b_{I} \\
\hline
\begin{smallmatrix}& & \circ\\
\circ &\circ &\circ &\circ &\circ &\circ &\circ &\\
\end{smallmatrix}
&25080&256& 
\begin{smallmatrix}& & \circ\\
\bullet &\circ &\circ &\circ &\circ &\circ &\circ &\\
\end{smallmatrix}
&4554&143\\
\hline
\begin{smallmatrix}& & \bullet\\
\circ &\circ &\circ &\circ &\circ &\circ &\circ &\\
\end{smallmatrix}
&6188&124& 
\begin{smallmatrix}& & \circ\\
\circ &\bullet &\circ &\circ &\circ &\circ &\circ &\\
\end{smallmatrix}
&7397&153\\
\hline
\begin{smallmatrix}& & \circ\\
\circ &\circ &\bullet &\circ &\circ &\circ &\circ &\\
\end{smallmatrix}
&9541&163& 
\begin{smallmatrix}& & \circ\\
\circ &\circ &\circ &\bullet &\circ &\circ &\circ &\\
\end{smallmatrix}
&8691&153\\
\hline
\begin{smallmatrix}& & \circ\\
\circ &\circ &\circ &\circ &\bullet &\circ &\circ &\\
\end{smallmatrix}
&7169&149& 
\begin{smallmatrix}& & \circ\\
\circ &\circ &\circ &\circ &\circ &\bullet &\circ &\\
\end{smallmatrix}
&5512&111\\
\hline
\begin{smallmatrix}& & \circ\\
\circ &\circ &\circ &\circ &\circ &\circ &\bullet &\\
\end{smallmatrix}
&4452&121& 
\begin{smallmatrix}& & \bullet\\
\bullet &\circ &\circ &\circ &\circ &\circ &\circ &\\
\end{smallmatrix}
&1422&58\\
\hline
\begin{smallmatrix}& & \circ\\
\bullet &\bullet &\circ &\circ &\circ &\circ &\circ &\\
\end{smallmatrix}
&1121&79& 
\begin{smallmatrix}& & \circ\\
\bullet &\circ &\bullet &\circ &\circ &\circ &\circ &\\
\end{smallmatrix}
&1981&85\\
\hline
\begin{smallmatrix}& & \circ\\
\bullet &\circ &\circ &\bullet &\circ &\circ &\circ &\\
\end{smallmatrix}
&1818&80& 
\begin{smallmatrix}& & \circ\\
\bullet &\circ &\circ &\circ &\bullet &\circ &\circ &\\
\end{smallmatrix}
&1554&85\\
\hline
\end{array}
$$

\newpage

\begin{center}
Table $E_{8}$ (cont'd)
\end{center}

$$
\begin{array}{||c|c|c||c|c|c||}
\hline
I & \sharp\mathcal{N}_{I} & \sharp \mathcal{A}b_{I} &
I & \sharp\mathcal{N}_{I} & \sharp \mathcal{A}b_{I} \\
\hline
\begin{smallmatrix}& & \bullet\\
\circ &\circ &\bullet &\circ &\bullet &\circ &\circ &\\
\end{smallmatrix}
&689&49& 
\begin{smallmatrix}& & \bullet\\
\circ &\circ &\bullet &\circ &\circ &\bullet &\circ &\\
\end{smallmatrix}
&580&39\\
\hline
\begin{smallmatrix}& & \bullet\\
\circ &\circ &\bullet &\circ &\circ &\circ &\bullet &\\
\end{smallmatrix}
&388&45& 
\begin{smallmatrix}& & \bullet\\
\circ &\circ &\circ &\bullet &\bullet &\circ &\circ &\\
\end{smallmatrix}
&553&49\\
\hline
\begin{smallmatrix}& & \bullet\\
\circ &\circ &\circ &\bullet &\circ &\bullet &\circ &\\
\end{smallmatrix}
&681&39& 
\begin{smallmatrix}& & \bullet\\
\circ &\circ &\circ &\bullet &\circ &\circ &\bullet &\\
\end{smallmatrix}
&475&49\\
\hline
\begin{smallmatrix}& & \bullet\\
\circ &\circ &\circ &\circ &\bullet &\bullet &\circ &\\
\end{smallmatrix}
&412&34& 
\begin{smallmatrix}& & \bullet\\
\circ &\circ &\circ &\circ &\bullet &\circ &\bullet &\\
\end{smallmatrix}
&412&44\\
\hline
\begin{smallmatrix}& & \bullet\\
\circ &\circ &\circ &\circ &\circ &\bullet &\bullet &\\
\end{smallmatrix}
&260&34& 
\begin{smallmatrix}& & \circ\\
\circ &\bullet &\bullet &\bullet &\circ &\circ &\circ &\\
\end{smallmatrix}
&566&43\\
\hline
\begin{smallmatrix}& & \circ\\
\circ &\bullet &\bullet &\circ &\bullet &\circ &\circ &\\
\end{smallmatrix}
&815&55& 
\begin{smallmatrix}& & \circ\\
\circ &\bullet &\bullet &\circ &\circ &\bullet &\circ &\\
\end{smallmatrix}
&668&47\\
\hline
\begin{smallmatrix}& & \circ\\
\circ &\bullet &\bullet &\circ &\circ &\circ &\bullet &\\
\end{smallmatrix}
&463&50& 
\begin{smallmatrix}& & \circ\\
\circ &\bullet &\circ &\bullet &\bullet &\circ &\circ &\\
\end{smallmatrix}
&726&45\\
\hline
\end{array}
$$

\newpage

\begin{center}
Table $E_{8}$ (cont'd)
\end{center}

$$
\begin{array}{||c|c|c||c|c|c||}
\hline
I & \sharp\mathcal{N}_{I} & \sharp \mathcal{A}b_{I} &
I & \sharp\mathcal{N}_{I} & \sharp \mathcal{A}b_{I} \\
\hline
\begin{smallmatrix}& & \bullet\\
\circ &\bullet &\bullet &\bullet &\circ &\circ &\circ &\\
\end{smallmatrix}
&105&16& 
\begin{smallmatrix}& & \bullet\\
\circ &\bullet &\bullet &\circ &\bullet &\circ &\circ &\\
\end{smallmatrix}
&187&22\\
\hline
\begin{smallmatrix}& & \bullet\\
\circ &\bullet &\bullet &\circ &\circ &\bullet &\circ &\\
\end{smallmatrix}
&161&18& 
\begin{smallmatrix}& & \bullet\\
\circ &\bullet &\bullet &\circ &\circ &\circ &\bullet &\\
\end{smallmatrix}
&112&21\\
\hline
\begin{smallmatrix}& & \bullet\\
\circ &\bullet &\circ &\bullet &\bullet &\circ &\circ &\\
\end{smallmatrix}
&255&28& 
\begin{smallmatrix}& & \bullet\\
\circ &\bullet &\circ &\bullet &\circ &\bullet &\circ &\\
\end{smallmatrix}
&301&24\\
\hline
\begin{smallmatrix}& & \bullet\\
\circ &\bullet &\circ &\bullet &\circ &\circ &\bullet &\\
\end{smallmatrix}
&211&30& 
\begin{smallmatrix}& & \bullet\\
\circ &\bullet &\circ &\circ &\bullet &\bullet &\circ &\\
\end{smallmatrix}
&192&22\\
\hline
\begin{smallmatrix}& & \bullet\\
\circ &\bullet &\circ &\circ &\bullet &\circ &\bullet &\\
\end{smallmatrix}
&186&28& 
\begin{smallmatrix}& & \bullet\\
\circ &\bullet &\circ &\circ &\circ &\bullet &\bullet &\\
\end{smallmatrix}
&116&22\\
\hline
\begin{smallmatrix}& & \bullet\\
\circ &\circ &\bullet &\bullet &\bullet &\circ &\circ &\\
\end{smallmatrix}
&121&23& 
\begin{smallmatrix}& & \bullet\\
\circ &\circ &\bullet &\bullet &\circ &\bullet &\circ &\\
\end{smallmatrix}
&184&19\\
\hline
\begin{smallmatrix}& & \bullet\\
\circ &\circ &\bullet &\bullet &\circ &\circ &\bullet &\\
\end{smallmatrix}
&124&25& 
\begin{smallmatrix}& & \bullet\\
\circ &\circ &\bullet &\circ &\bullet &\bullet &\circ &\\
\end{smallmatrix}
&165&23\\
\hline
\end{array}
$$

\newpage

\begin{center}
Table $E_{8}$ (cont'd)
\end{center} 

$$
\begin{array}{||c|c|c||c|c|c||}
\hline
I & \sharp\mathcal{N}_{I} & \sharp \mathcal{A}b_{I} &
I & \sharp\mathcal{N}_{I} & \sharp \mathcal{A}b_{I} \\
\hline
\begin{smallmatrix}& & \circ\\
\bullet &\bullet &\circ &\circ &\bullet &\bullet &\bullet &\\
\end{smallmatrix}
&27&10& 
\begin{smallmatrix}& & \circ\\
\bullet &\circ &\bullet &\bullet &\bullet &\bullet &\circ &\\
\end{smallmatrix}
&48&11\\
\hline
\begin{smallmatrix}& & \circ\\
\bullet &\circ &\bullet &\bullet &\bullet &\circ &\bullet &\\
\end{smallmatrix}
&54&15& 
\begin{smallmatrix}& & \circ\\
\bullet &\circ &\bullet &\bullet &\circ &\bullet &\bullet &\\
\end{smallmatrix}
&47&15\\
\hline
\begin{smallmatrix}& & \circ\\
\bullet &\circ &\bullet &\circ &\bullet &\bullet &\bullet &\\
\end{smallmatrix}
&36&13& 
\begin{smallmatrix}& & \circ\\
\bullet &\circ &\circ &\bullet &\bullet &\bullet &\bullet &\\
\end{smallmatrix}
&26&10\\
\hline
\begin{smallmatrix}& & \bullet\\
\circ &\bullet &\bullet &\bullet &\bullet &\circ &\circ &\\
\end{smallmatrix}
&32&10& 
\begin{smallmatrix}& & \bullet\\
\circ &\bullet &\bullet &\bullet &\circ &\bullet &\circ &\\
\end{smallmatrix}
&49&9\\
\hline
\begin{smallmatrix}& & \bullet\\
\circ &\bullet &\bullet &\bullet &\circ &\circ &\bullet &\\
\end{smallmatrix}
&35&12& 
\begin{smallmatrix}& & \bullet\\
\circ &\bullet &\bullet &\circ &\bullet &\bullet &\circ &\\
\end{smallmatrix}
&58&12\\
\hline
\begin{smallmatrix}& & \bullet\\
\circ &\bullet &\bullet &\circ &\bullet &\circ &\bullet &\\
\end{smallmatrix}
&58&15& 
\begin{smallmatrix}& & \bullet\\
\circ &\bullet &\bullet &\circ &\circ &\bullet &\bullet &\\
\end{smallmatrix}
&36&12\\
\hline
\begin{smallmatrix}& & \bullet\\
\circ &\bullet &\circ &\bullet &\bullet &\bullet &\circ &\\
\end{smallmatrix}
&68&12& 
\begin{smallmatrix}& & \bullet\\
\circ &\bullet &\circ &\bullet &\bullet &\circ &\bullet &\\
\end{smallmatrix}
&75&18\\
\hline
\end{array}
$$

\newpage

%% Table F4

\begin{center}
Table $F_{4}$
\end{center}

$$
\begin{array}{||c|c|c||c|c|c||}
\hline
I & \sharp\mathcal{N}_{I} & \sharp \mathcal{A}b_{I}  & 
I &  \sharp\mathcal{N}_{I} &  \sharp\mathcal{A}b_{I} \\
\hline
\begin{smallmatrix}& & 
\circ &\circ & \circ &\circ & \\
\end{smallmatrix}
&105&16& 
\begin{smallmatrix}& & \\
\bullet &\circ &\circ &\circ &\\
\end{smallmatrix}
&24&6\\
\hline
\begin{smallmatrix}& & 
\circ &\bullet & \circ &\circ &\\
\end{smallmatrix}
&35&12& 
\begin{smallmatrix}& & \\
\circ &\circ & \bullet &\circ & \\
\end{smallmatrix}
&32&10 \\
\hline
\begin{smallmatrix}& & 
\circ &\circ &\circ &\bullet &\\
\end{smallmatrix}
&49&9&
\begin{smallmatrix}& & \\
\bullet &\bullet & \circ &\circ & \\
\end{smallmatrix}
&10&5 \\
\hline
\begin{smallmatrix}& & 
\bullet &\circ &\bullet &\circ &\\
\end{smallmatrix}
&8&4&
\begin{smallmatrix}& & \\
\bullet &\circ & \circ &\bullet &\\
\end{smallmatrix}
&12&4 \\
\hline
\begin{smallmatrix}& & 
\circ &\bullet &\bullet &\circ &\\
\end{smallmatrix}
&14&7 &
\begin{smallmatrix}& & \\
\circ &\bullet & \circ &\bullet & \\
\end{smallmatrix}
&14&6 \\ 
\hline 
\begin{smallmatrix}& & 
\circ &\circ &\bullet &\bullet &\\
\end{smallmatrix}
&10&4&
\begin{smallmatrix}& & \\
\bullet &\bullet & \bullet &\circ & \\
\end{smallmatrix}
&4&3 \\
\hline
\begin{smallmatrix}& & 
\bullet &\bullet &\circ &\bullet &\\
\end{smallmatrix}
&5&3 &
\begin{smallmatrix}& & \\
\bullet &\circ & \bullet &\bullet &\\
\end{smallmatrix}
&3&2 \\
\hline
\begin{smallmatrix}& & 
\circ &\bullet &\bullet &\bullet &\\
\end{smallmatrix}
&3&2&
\begin{smallmatrix}& & \\
\bullet &\bullet & \bullet &\bullet & \\
\end{smallmatrix}
&1&1 \\
\hline
\end{array}
$$

where we use the following orientation for the Dynkin diagram of $F_4$ :

$$
\begin{Dynkin}
\Dbloc{\Dcirc\Deast\Dtext{t}{1}}
\Dbloc{\Dcirc\Dwest\Ddoubleeast\Dtext{t}{2}}
\Drightarrow
\Dbloc{\Dcirc\Ddoublewest\Deast\Dtext{t}{3}}
\Dbloc{\Dcirc\Dwest\Dtext{t}{4}}
\end{Dynkin}
$$

%% Table G2

\begin{center}
Table $G_{2}$
\end{center}

$$
\begin{array}{||c|c|c||c|c|c||}
\hline
I & \sharp\mathcal{N}_{I} & \sharp \mathcal{A}b_{I}  & 
I &  \sharp\mathcal{N}_{I} &  \sharp\mathcal{A}b_{I} \\
\hline
\begin{smallmatrix}& & 
\circ &\circ & \\
\end{smallmatrix}
&8&4& 
\begin{smallmatrix}& & \\
\bullet &\circ &\\
\end{smallmatrix}
&3&2\\
\hline
\begin{smallmatrix}& & 
\circ &\bullet & \\
\end{smallmatrix}
&4&3& 
\begin{smallmatrix}& & \\
\bullet &\bullet &\\
\end{smallmatrix}
&1&1\\
\hline
\end{array}
$$

where we use the following orientation for the Dynkin diagram of $G_2$ :

$$
\begin{Dynkin}
\Dbloc{\Dcirc\Deast\Ddoubleeast\Dtext{t}{1}}
\Dleftarrow
\Dbloc{\Dcirc\Dwest\Ddoublewest\Dtext{t}{2}}
\end{Dynkin}
$$

\end{document}